\documentclass[10pt,journal]{IEEEtran}

\usepackage{graphicx}
\usepackage{dcolumn}
\usepackage{bm}
\usepackage{hyperref}
\usepackage[mathlines]{lineno}
\usepackage{dsfont} 
\usepackage{multirow}
\usepackage{amsmath}    
\usepackage{amsfonts}  

\newenvironment{definition}[1][Definition]{\begin{trivlist}
\item[\hskip \labelsep {\bfseries #1}]}{\end{trivlist}}

\begin{document}


\title{Chaotic Encryption Method Based on Life-Like Cellular Automata}
\author{
\IEEEauthorblockN{Marina Jeaneth Machicao, Anderson G. Marco, Odemir M. Bruno\\ }
    \IEEEauthorblockA{Instituto de F\'{i}sica de S\~{a}o Carlos - Universidade de S\~{a}o Paulo - Brasil
    \\mj.machicao@gmail.com ; anderson.marco@gmail.com ; bruno@ifsc.usp.br}
}



\maketitle

\begin{abstract}
We propose a chaotic encryption method based on Cellular Automata(CA), specifically on the family called the ``Life-Like'' type. Thus, the encryption process lying on the pseudo-random numbers generated (PRNG) by each CA's evolution, which transforms the password as the initial conditions to encrypt messages.
Moreover, is explored the dynamical behavior of CA to reach a ``good'' quality as PRNG based on measures to quantify ``how chaotic a dynamical system is'', through the combination of the entropy, Lyapunov exponent, and Hamming distance. Finally, we present the detailed security analysis based on experimental tests: DIEHARD and ENT suites, as well as Fourier's Power Spectrum, used as a security criteria.
\end{abstract}

\section{Introduction}
Cryptography can be traced back over 2300 years as decisive in the course of human History.
The sort of battles and kings were decided by the power or weakness of ciphers.
Nowadays, its importance lies on: military communications, intelligence tool, banking system, and e-commerce foundation.
However in the last 20 years many classic encryption algorithms had been broken, among
them, DES algorithm in 1998\cite{CrackDES}, cryptography hash functions MD5
(2008)\cite{CrackMD5} and SHA1 (2009)\cite{CrackSHA1}.
Furthermore, Cryptography is constantly searching for new algorithms to encrypt and to make ciphers more secure, some new mathematical methodologies has to be
considered.

Chaos Theory is a mathematical field that studies the behavior of complex systems with highly sensitiveness to initial conditions. In fact, chaos can be found in very simple systems that shows complex patterns, for instances, Chaos was first introduce in 1960 as a weather prediction by Edward Lorenz, furthermore  it was involved with population growth, and also with cellular evolution introduced as Cellular Automata. In sum Chaos Theory deals with unpredictable systems, most of them with systems found in nature: weather, turbulence, stock market, population, and so on. 

Cellular automata, first introduced by von Neumann in the early 1950s, was conceived as a
model of biological evolution \cite{neumann}, and also used as a prototyping model for a large variety of natural systems, this models have attracted the attention of several groups of research as CA provide approximation to partial differential equations, and because CA are considered discrete dynamical systems and furthermore computational systems.
Twenty years later, the mathematician Jhon Conway introduced his well-known Life game, which initially was intended to modeling natural cellular phenomena, consists of a regular grid of cells, with alive-dead states(``On'' and ``Off''), through out simple rules, which involves neighborhood`s cells, they constantly evolved even to get unpredictable behavior\cite{Aguena, Wolfram}.

Over the past decade, there exists progress on research encryption techniques based on Chaos Theory\cite{Xiang, Anderson, Deng}. This position was motivated by the similarity between Chaos and Cryptography, where were emphasized the sensitiveness of initial conditions trough the \emph{key} and input(plaintext) dependency, also both Chaotic-System and Crypto-System has a ``random-like'' behavior\cite{Alvarez2006}, which is very desirable for cryptographic schemes. By extension all this relationship are inherited to Cellular Automata. Consequently cryptography can take advantages from those CA which shows chaotic behavior.

Cryptography based on CA gave first signals at \cite{Wolfram, Habutsu, Nandi} and recently dealing in \cite{Tomassini,Seredynski}, it was also included CA as pseudo-random number generator in \cite{TomassiniRandom, Martin,Guan, Shin, Shin2009, Seth}.
Although the inclusion of CA in the field of cryptography is not new, this proposal follows a chaotic sensitivity throughout the process, and our results exceed expectations compared to others found in literature.

This motivated us to develop an chaos-based encryption scheme based on CA that takes advantage of chaos to encrypt in a ``unpredictable'' manner, this provides a very attractive alternative method as there exists CA rules with chaotic behavior that, can be employed as pseudo-random number generators, and because it shows efficiently hardware incorporation according to its dynamical properties and through its parallelism fundamentals. To understand what this ``chaotic behavior'' really means, where used the Lyapunov exponent as a measure that characterizes the rate of separation of infinitesimally close trajectories provided by the system, in advance it determines a notion of predictability for a dynamical system, where positives values is usually taken as a chaotic indicator. \cite{Wolf} 

In Section 2, we mention the state of art of some drawbacks based on chaos-based cryptography. Section 3 introduces basic theory on the subject, as needed. In Section 4, the encryption method based on CA is proposed in detail.

Section 5, presents the security analysis based on experimental tests: DIEHARD\cite{DIEHARD} and ENT\cite{ENT}, as well as Fourier's Power Spectrum, used as a security criteria. Finally, in Section 6, the paper ends with a discussion about the method.

\section{Background}

Modern cryptography is based on two approaches: Discrete Mathematics (Symmetric key) and Number Theory (Asymmetric key)\cite{Tanenbaum, Stallings}. Cryptography's algorithms stand on an elementary mathematics need to be extremely complex (Symmetric key Algorithm), performing many bit-wise operations and permutation between neighbors, on the other hand, having long keys (Asymmetric key Algorithm). However, for both strategies, computing complexity and the size of the keys, just make information to be secure for a certain time.

Since Baptista has introduced first chaos-based encryption\cite{Baptista}, many researches have been interested in the relationship between chaos and cryptography, then many properties of chaotic systems have to be on mind e.g: ergodicity, sensitivity to
initial conditions, deterministic dynamics, and structural complexity.

Cryptography of symmetric key systems based on CA,  were first studied by Wolfram \cite{Wolfram}, Habutsu\cite{Habutsu}, Nandi et al.\cite{Nandi} and Gutowitz\cite{Gutowitz}, and later by Tomassini et al.\cite{Tomassini}. Recently one-dimensional CA was subject of study by Seredynski et al.\cite{Seredynski}. This paper  presents a new effort to involve Chaos Theory and Cellular Automata to implements an encryption method.
\section{Cellular Automata}
\label{SECTION:CA}
In this work a cellular automata (CA) is considered a discrete dynamical system
defined on a discrete space which is governed by its local rules and by its immediate neighbors, which specifies how CA evolves in time.

\begin{definition}
An homogeneous CA can be represented as a sextuple $\langle \mathds{T},S,s, s_{0},N,\phi \rangle$ where:
\begin{enumerate}
  \item $\mathds{T}$ is a lattice of a n-dimensional Euclidean space $\mathds{R}^n$, consisting of cells $c_{i}$, $i \in \mathds{N}$
  \item S is a finite set of $k$ states, often $S \subset \mathds{N}$
  \item The output function $s: \mathds{T}\times \mathds{N} \rightarrow S$ maps the states of cell $c_{i}$ at discrete time $t$, i.e., $s(c_{i},t)$.
  \item The function $s_{0}: A \rightarrow S$ allocates the initial configuration of every cell $c_{i}$, i.e., $s(c_{i})= s_{0}(c_{i})$.
  \item The neighborhood function $N: \mathds{T} \rightarrow \bigcup_{p=1}^{\infty}\mathds{T}^{p}$, yields every cell $c_{i}$ to a finite sequence $N(c_{i})=(c_{i_{j}})_{j=1}^{|N(c_{i})|}$, with $|N(c_{i})|$ distinct cells $c_{i}$.
  \item The transition function $\phi:S^{|N(c_{i})|}\rightarrow S$, describes the rules governing the dynamics of every cell $c_{i}$
\begin{equation}s(c_{i}, t+1) = \phi((s(c_{i},t))_{j=1}^{|N(c_{i})|}) = \phi(\sigma_{i}) \nonumber \end{equation} with,
\begin{equation} \sigma_{i} = \sum_{j=1}^{N(c_{i})}s(c_{i},t)  \end{equation}
\end{enumerate}
\end{definition}

\subsection{Life-Like Cellular Automata}
A well known two-dimensional CA was proposed by John Conway, called ``Game of Life'', which is based on a biological model, with cells either dead$(0)$ or alive$(1)$, in general the state of a cell at the next generation depends on its own state and the sum of cells, to Conway's CA a dead cell births when is surrounding for 3 alive neighbors and an alive cell can survive if there exists 2 or 3 active cells.

Using a standard convention of Golly simulator\cite{Golly}, a notation for naming CA, a rule is written in the form $By\backslash Sx$ where $x$ and $y$ is a sequence of distinct digits from 0 to 8, in numerical order. Thus, Conway's Game of Life is denoted $B3\backslash S23$, where ``B'' stands for ``birth'' and the ``S'' stands for ``survival''. It is common to refer to ``Life-Like'' or simply ``Life Family'' in the sense of similarity to Conway's, to all those CA in the format of Golly that have one of several rules, for instances $(B36\backslash S23)$, $(B23\backslash S36)$ .

Formally, any ``Life-Like'' automata is described as $\langle \mathds{Z}^{2},S=\{0,1\},s, s_{0},$ Moore neighborhood $,\phi: S^{9} \rightarrow S \rangle$, where Moore neighborhood considers 8 cardinal direction and the state of the own center. For instances, for binary cells $c_{1},c_{2},c_{3}, c_{4}, c_{5}, c_{6}, c_{7},c_{8}$, and $c_{9}$ we say that the transition function, at any time $t$, for rule $B3\backslash S23$(Game-of-Life)\cite{Adamatzky}  is of the form:

\begin{equation}
    \phi \Bigg(\begin{array}{l l l}
        c_{1}, & c_{2}, & c_{3},\\
        c_{4}, & c_{5}, & c_{6},\\
        c_{7}, & c_{8}, & c_{9}\\
      \end{array} \Bigg)=
    \left\{
      \begin{array}{l l}
        1, & \quad \text{if }         \displaystyle\sum_{i=1}^{9} s(c_{i},t) = 3\\
        1, & \quad \text{if }         \displaystyle\sum_{i=1}^{9} s(c_{i},t)= 2, i\neq5\\
        0, & \quad \text{otherwise}
      \end{array} \right.
\end{equation}

There are two main issues with this ``Life-Family'' CA:
\begin{itemize}
  \item \textbf{Cellular growth and decay}: The population of cells of some rules of the Life-Family seems to decay when iterated considerable time. It is mean, they are ``like-death'' because most of their cells are \emph{like} to arise death\cite{BookChaos}.
  \item \textbf{Boundary conditions}: By definition, a 2D-CA consists of an \emph{infinite} plane, for instances, its space $\mathds{T}$ can be a matrix of order $m\times n$. For computational reasons they are often simulated on a finite grid rather than an infinite one, although there is an evidently issue with the boundaries. There exists plenty methods to handled it, in this paper is used a hyper-toroidal method,  where the grid is considered as their edges were touching on all boundaries. So far, this toroidal arrangement simulates an infinite periodic lattice.
\end{itemize}

\subsection{Chaotic Cellular Automata}
In some cases, CAs present chaotic behavior, then we are interested to use chaotic CAs in order to take advantage for cryptography. To measure the chaotic behavior in dynamical systems were announced several studies that involved both geometrical\cite{Guckenheimer} and statistical\cite{EckmannRuelle} approaches. The statistical approach, seeks to characterize dynamical systems through the Lyapunov exponent, which has been proven to be the most useful to measure chaos\cite{Tisseur2003,Tisseur2005,Cattaneo,XUXU}.

\subsubsection{Lyapunov Exponent(LE)}
The Lyapunov exponent $\lambda$  is a measure of the sensitive dependence on initial conditions, it is a measure of how chaotic a dynamical system is.
There are several ways to estimate the LE, and several variants were suggested\cite{Shimada,Benettin,EckmannRuelle,Guckenheimer}. Where, the information extracted from  $\lambda$ is:
 \begin{equation}
\left\{
      \begin{array}{l l}
        \lambda < 0 & \quad \text{Stable Periodic}\\
        \lambda > 0 & \quad \text{Chaotic}\\
        \lambda = 0 & \quad \text{Neutrally Stable}
      \end{array} \right.
\end{equation}

For this purpose is considered a recent method to estimate the LE in CA, according with Baetens and De Baets \cite{Baetens2010}, where is considered  two initial configurations $s_{0}$ and $s^*_{0}$ of CA with states $s=\{0,1\}$. Is defined a ``damage vector'', $h(\cdot, t) = s(\cdot,t)\oplus s^*(\cdot,t)$, as the number of different cells or distance between both CA. Is also defined a \textbf{perturbed cell} $c_{i}$ for which $s(c_{i},t)\neq s^*(c_{i},t)$ and the smallest perturbation of $s_{0}$ due to the state domain of CA, $\sum_{c_{i}}h(c_{i},0)=1$.

Is is important to take into account all the perturbations originated from one initial perturbation, where can be defined as the maximum Lyapunov exponent $\lambda(t)$, with finite $t \in \mathds{N}$

\begin{equation}
\lambda(t) = \frac{1}{t}\log(\frac{\sum_{c_{i}}h(c_{i},t)}{\sum_{c_{i}}h(c_{i},0)})
\end{equation}

\subsubsection{Entropy}
Entropy is a measure of disorder or randomness in a closed system. Entropy applied to CA can be estimated by:

 \begin{equation}
 \begin{array}{l}
H = - \sum_{i=0}^{k}P_{i}\log(P_{i})
\\
\\
0\leq H\leq 1
\end{array}
\end{equation}

Where $P_{i}$ is the probability of $k$ possible states of every cell in CA. By using entropy, we are interested on maximal values, thus makes $1$ the maximum normalized entropy,  which would mean a low level of redundancy or predictability of the encryption method.

\subsubsection{Hamming Distance (HD)}
We employ the Hamming distance $D_{H}$ in CA as a measure that estimates the minimum number of substitutions required to change one state $c_{i}$ into another, this means the count of differences between the states of cells $s(\cdot,t)$ and $s(\cdot,t+1)$, then, as the higher disturbances the lower repetition of patterns.
 \begin{equation}
D_{H} = \frac{\sum s(\cdot,t)\oplus s(\cdot,t+1)} {size(\mathds{T})}
\end{equation}
Finally, is considered the highest Hamming distance, when all cells are the opposed to next iteration, thus the maximal $D_{H}$ is closed to $1$.

\subsubsection{Chaos Combined Measure (Max)}
It is notice that isolated high Lyapunov exponent, entropy or Hamming distance are necessary, but by no means enough to response the question: ``which CA is more chaotic than the other?''. Then, this new measure \emph{Max} maximize the combination of the measures mentioned above, by simple multiplicity. Moreover, as the higher \emph{Max} gets, the more chaotic a system is. this means:
 \begin{equation}
\text{Max} = (\lambda * H  *D_{H})
\end{equation}

\section{The proposed Algorithm}
We propose a chaotic encryption method based on Cellular Automata(CA), specifically on two-dimensional ``Life-like'' CA, mentioned on section III, to design a symmetric key cryptography system.

\begin{figure*}[hbtp]
\begin{center}
\includegraphics[width=1.0\textwidth,height=0.42\textwidth]{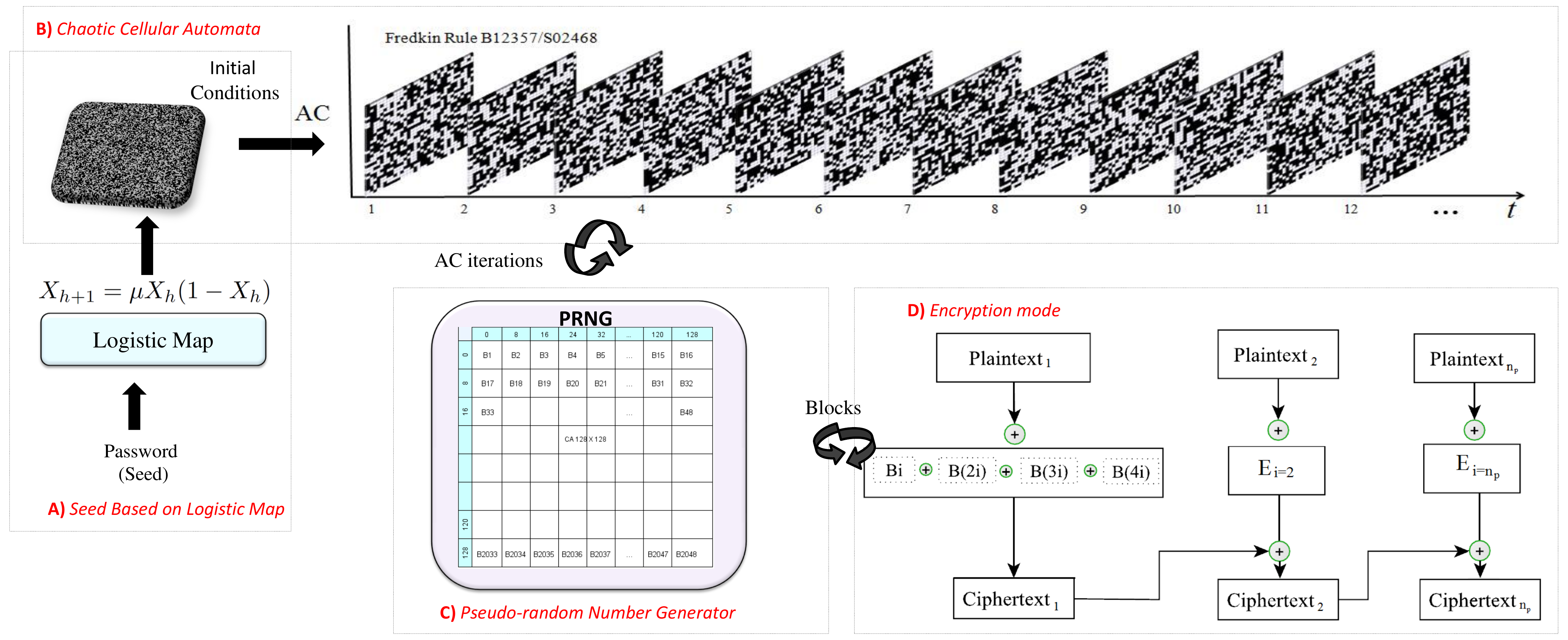}
\caption{Proposal of chaotic encryption method based on Cellular
Automata(CA), for instances, it uses 4 Blocks as PRNG. Its design is divided into 4 principal parts. (A) The transformation of the password into a initial CA, (B) The CA is defined with a specific rule and a initial CA and is iterated as many as required (C) Transforms the CA into a PRNG to generate blocks, (D)The encryption mode takes blocks as many as needed. Note that the decryption uses the same design just than inverse the encrypt process.}
\label{fig:Proposta}
\end{center}
\end{figure*}

In Figure \ref{fig:Proposta} is showed the proposal, which is divided into 4 principal parts, where a chaotic CA is applied as a pseudo-random number generator (PRNG) which is employ during the encryption. With this process, is constantly taken pseudo-random numbers generated at every step during the CA's evolution, this numbers are composed in blocks of size of the plaintext, for later encryption. The initial conditions for this CA can be fulfilled in different ways, for this proposal was taken the Logistic Map (dynamical system) where the password to encrypt is taken as a seed for the system, thus, every cell takes binaries values as the Logistic Map iterates.
At next, is showed each of the 4 the parts of our cryptosystem:

\begin{itemize}
  \item Seed based on Logistic Map
  \item Chaotic  Cellular Automata
  \item Pseudo-random number generator
  \item Encryption/Decryption
\end{itemize}

\subsection{Seed based on Logistic Map}
The Logistic Map is a well-known continuous dynamical systems, mathematically is written as
\begin{equation}
\label{eq:Logistic}
X_{h+1} = \mu X_{h}(1-X_{h})
\end{equation}

Where $X_{h} \in [0,1]$ for $h$ $\in$ $\mathbf{N}$, with parameter $\mu$, which presents chaotic behavior for values among $3.9\leq \mu \leq 4.0 $.

This Logistic Map generates continuous values between $[0,1]$, which are discretized(binarized) in order to fulfill the initial cellular automata to later encryption.

This discretization is allowed because of the equally distribution of the Logistic Map
, this property is very important in order to be used in cryptography, thus, the probabilities of any password would also fit in an equally distributed region.

As described above, this map takes the password as a \textbf{seed} to generate discrete $X_{h}$, through it iterates until all cells in CA gets completed. Should be noted that have to be omitted first iterations in order to get rid of the transient  $\alpha$ of the map, which presents non interested region of Logistic Map.

To transform $X_{h}$ into discrete values, we have to consider the following definitions:

\begin{definition}
Let $\vec{\pi}$ be a vector containing a string(password), in bytes, with $size(\vec{\pi}) := 16$, this means length of 128 bits.
Is defined:
\begin{equation}
\Omega:= \displaystyle\sum_{i=1}^{size(\vec{\pi})} 2^{8(i-1)}\pi_{i}
\end{equation}
\begin{equation}
\Omega:= \frac{\Omega}{2^{(8m+1)}}.
\end{equation}
\end{definition}

Where $X_{0}= \Omega +\epsilon$, with $\Omega\in [0,1]$ and for very small $\epsilon$, as Logistic Map cannot take zero as an entrance value.

\begin{definition}
Let the space of CA  $\mathds{T}$ be a matrix of order $m\times n$, with entries $1\leq i \leq m; 1\leq j \leq n$, which also means $size(\mathds{T}) = m\times n$.

Let $\alpha$ $\in$ $\mathbf{N}$ be the transient of the Logistic map.
\end{definition}

\begin{equation}
  s_{0}(\cdot) = \left\{
  \begin{array}{l l}
    1 & \quad \text{if $X_{((m(i-1)+j) + \alpha)}  < 0.5$ }\\
    0 & \quad \text{if $X_{((m(i-1)+j) + \alpha)}  \geq 0.5$ }\\
  \end{array} \right.
\end{equation}

\subsection{Chaotic Cellular Automata}
In this section is discussed the different types of rules of the ``Life-Like'' family, as several rules can be considered to cryptography. In Table \ref{Tab:Fredkin} can be visually appreciated the random nature of 64x64 2D-CA ``Fredkin'' with rule $B1357\backslash S02468$, each of the 3 first figures present different percentage of initial seed of $10\%$, $50\%$, and $90\%$ of alive cells. Curiously, as seen, this rule stabilizes to $50\%$ of alive cells after few iterations and keep it constantly through time, even with seeds of $10\%$ or $90\%$. Furthermore, the measure \emph{Max} (combination of entropy, Lyapunov Exponent, and Hamming distance) mentioned in Section III, is then considered as a way to choose the most appropriate chaotic CA, which according to results showed in next section, the rule ``Fredkin'' is considered the most acceptable CA to cryptography.
\begin{table}[htbp]
\caption{2D-CA ``Fredkin''$B1357\backslash S02468$ evolution with different seeds}
\begin{tabular}{c}
\includegraphics[width=0.46\textwidth,height=0.17\textwidth]{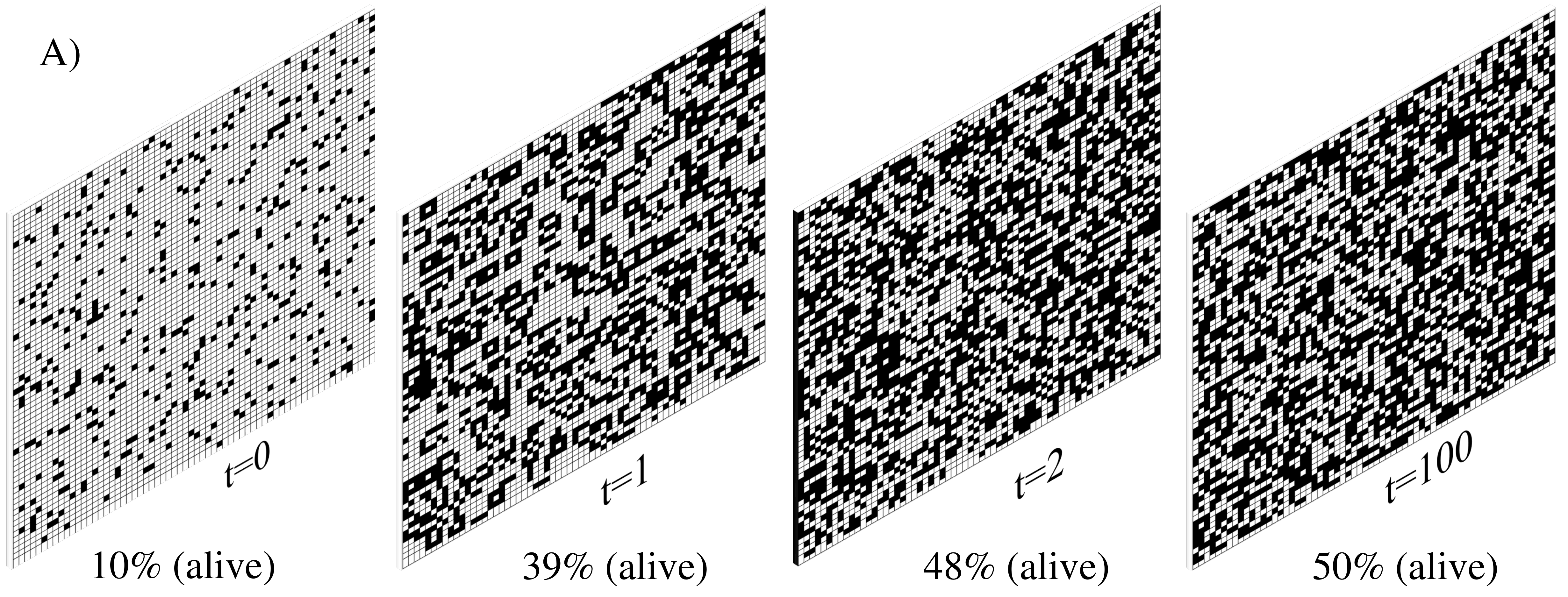}
\\
\hline
\\
\includegraphics[width=0.46\textwidth,height=0.17\textwidth]{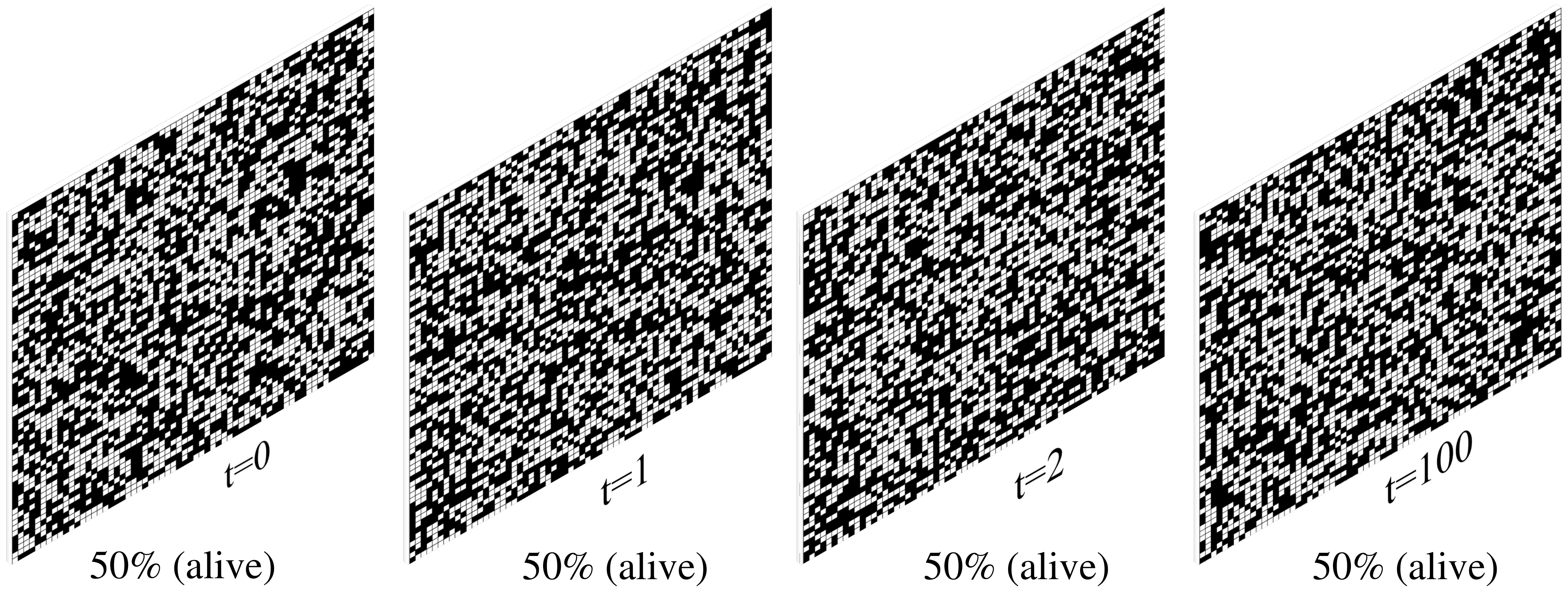}
\\
\hline
\\
\includegraphics[width=0.46\textwidth,height=0.17\textwidth]{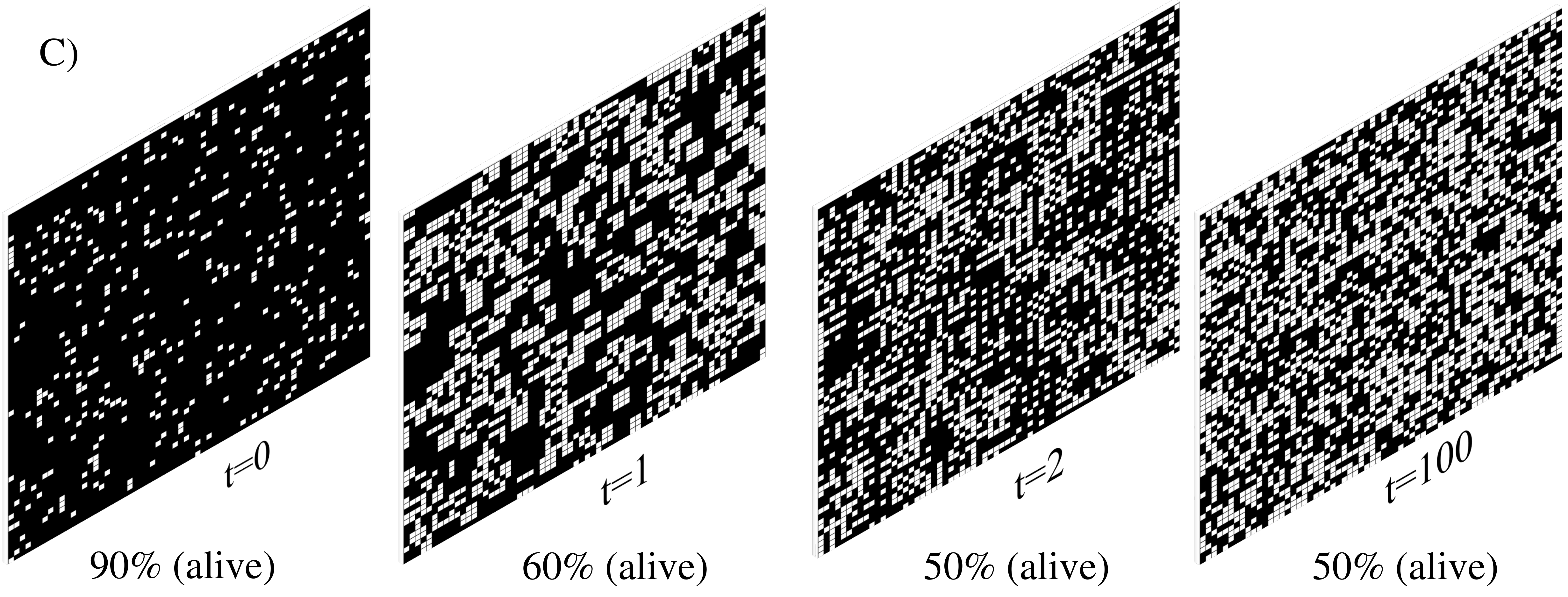}
\\
{\scriptsize 2D-CA Fredkin with random seed A)10\%, B)50\% and C)90\% of alive cells}
\label{Tab:Fredkin}
\end{tabular}
\end{table}

\subsection{Pseudo-random number generator(PRNG)}

A pseudo-random number generator is an algorithm for generating sequences of numbers that approximates randomness, this property  is required in different problems  domains such as simulations, testing, games, security, and of course making it suitable for use in cryptography. This approach is not recent, and exists several implementations in literature, the most recents ones are mentioned here in \cite{TomassiniRandom, Martin,Guan, Shin2009, Seth}, for this, a CA can be taken in several different ways.
The states of cells $c_{i}$ of an $m\times n$ CA  are binary values, then it can be taken as an advantage, by \textbf{composing blocks} which is described at next:
\begin{itemize}
\item Let $\vec{\beta}$ be a vector of bits (0,1)
    \begin{itemize}
        \item $size(\vec{\beta}) = floor(\frac{size(\vec{\beta})}{8})$
        \begin{equation}
        \beta_{(i-1)+((j-1)\times n)+ ((t-1)\times m\times n)+1} = s(\cdot,t)
        \end{equation}
    \end{itemize}

\item Let $\vec{B}$ be a vector of blocks of bytes

    \begin{equation}
    B_{i}= \displaystyle\sum _{j=1}^{8} 2^{8(j-1)}\beta_{8(i-1)+j}
    \end{equation}

\item Let $\vec{Y}$ be a vector composed by $\rho$ blocks.
    \begin{itemize}
        \item $size(\vec{Y}) \leq floor(\frac{size(\vec{B})}{8})$
    \end{itemize}
\end{itemize}

This ``composing-block'' process is defined as the consecutive XORed of the blocks $\vec{B}$ as follows:
\begin{equation}
Y_{i}= B_{(i\times \rho) +1}\oplus B_{(i\times \rho)+2} \oplus \ldots \oplus B_{(i\times \rho)+\rho}
\end{equation}

The symbol $\oplus$  represents the bitwise XOR logical operation, which is due to its reversible property.

\subsection{Encryption and Decryption}
Let $\overrightarrow{P}$ be the plaintext message of size $n_{p}$ and $E$ our ciphering algorithm, which represents any symmetrical cryptography algorithm using password $\vec{\pi}$. The fundamental transformation to obtain the ciphertext  $\overrightarrow{C}$ of same size $n_{p}$, is by encrypting the combined blocks generated by the PRNG and then XORed with  the previous plaintext sequence.

This procedural can ensure that similar blocks of the plaintext are codified as different ciphertext, making the cipher stronger.

For the Encryption mode:
\begin{equation}
\begin{array}{l l}
    C_{i}= E_{i}(P_{i} \otimes C_{i-1} \otimes Y_{i} ) & \quad \\
    C_{0}= Y_{1} \otimes P_{1} & \quad \text{ for } i=1,2,\ldots,n_{p}
\end{array}
\end{equation}

For the Decryption mode:
\begin{equation}
\begin{array}{l l}
    P_{i}= E_{i}^{-1}(C_{i}\otimes  Y_{i}) \otimes C_{i-1} & \quad  \\
    C_{0}= Y_{1} \otimes P_{1} & \quad  \text{ for } i=1,2,\ldots,n_{p}

\end{array}
\end{equation}


\section{Experimental results}

In this section, the efficiency of the chaotic encryption method based on Life-Like
CA is analyzed, thus it was separated into parts to presents the results:

\subsection{Chaotic Life-Like CA}
The purpose of this experiment was to discover among an enlarged set of rules of CAs, which one
presents the most chaotic behavior by having high entropy, Lyapunov Exponent, and Hamming distance,
in order to produce very high quality of PRNG and furthermore a ``good'' cryptosystem.

In Table \ref{Tab:LifeLike} are showed the results of the test with a  2D-CA of  128x128 cells of
the Life-Like family, to estimate the average entropy were employed 20 million iterations, the
LE was estimated under $T=200$, and the average Hamming distance was calculated with $T=1000$.
Besides, the maximal values were taken as the highest combination $Max$, this means, by
selecting the ``Fredkin'', ``Amoeba'', and ``$B23\backslash S36$'' rules as some of the best
chaotic CA to be employed for this encryption proposal.

\begin{table}[htbp]
\caption{Comparison of Life-Like rules, based on its entropy, Lyapunov Exponent, and Hamming distance}
\begin{tabular}{c}

\includegraphics[width=0.5\textwidth,height=0.48\textwidth]{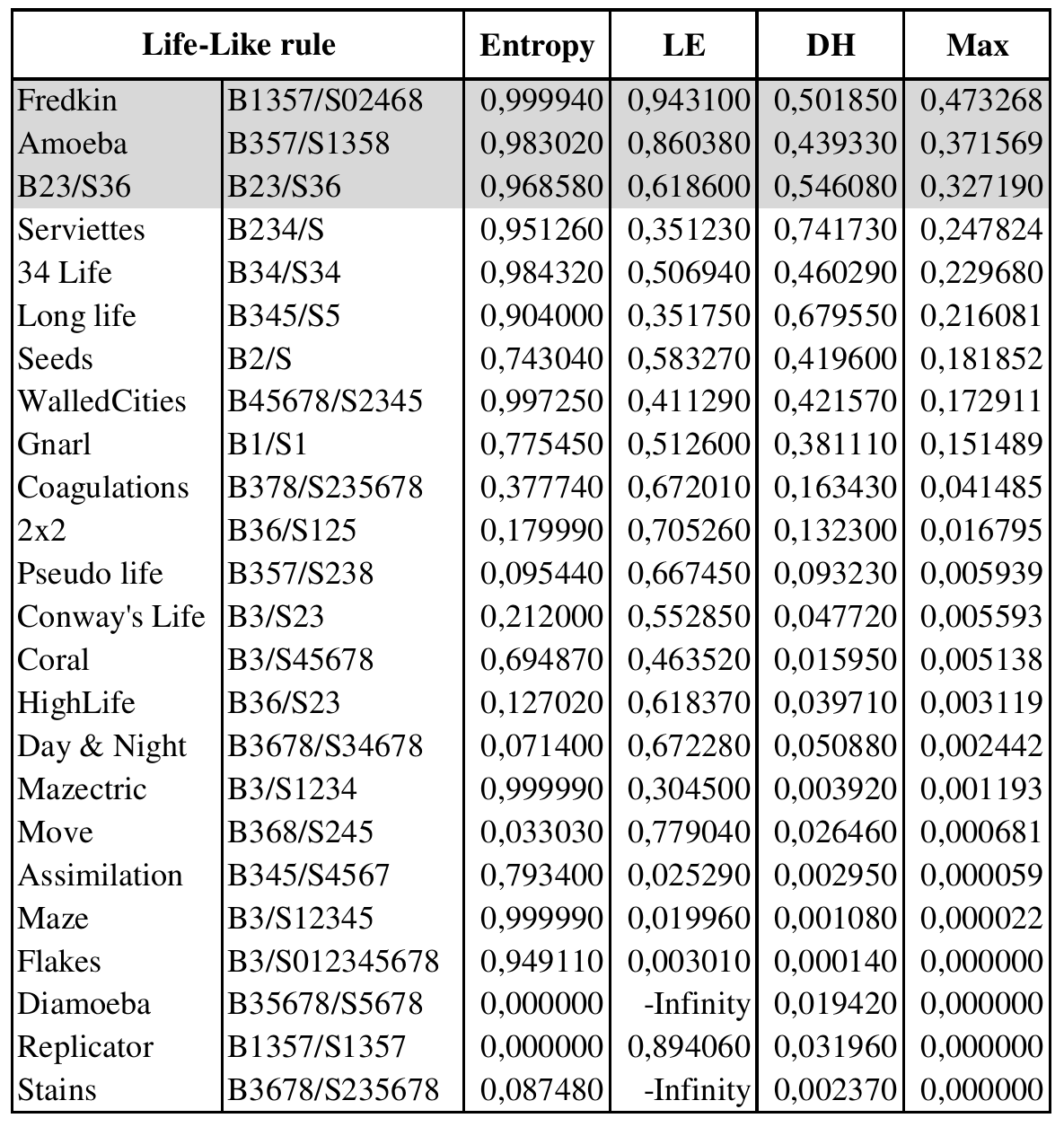}

\label{Tab:LifeLike}
\end{tabular}
\end{table}

For all subsequent experiment in this paper, is used the rule \emph{Fredkin} $B1357\backslash S02468$, as it posses constantly and highly entropy, that also presents chaotic behavior, and finally guarantee a higher disturbance variation as it is iterated, then this model is applied in this proposal. the random nature of this \emph{Fredkin} rule can be visually appreciated in Table \ref{Tab:LifeLikeComplete}, where is presented the complete analyze of the chaotic measures presents above, the pictures of each CA, were take at $t=1000$.

\begin{table*}[htbp]\begin{center}
\caption{Life-Like rules, based on its entropy, Lyapunov Exponent, and Hamming distance}
\begin{tabular}{c}

\includegraphics[width=0.9\textwidth,height=1\textwidth]{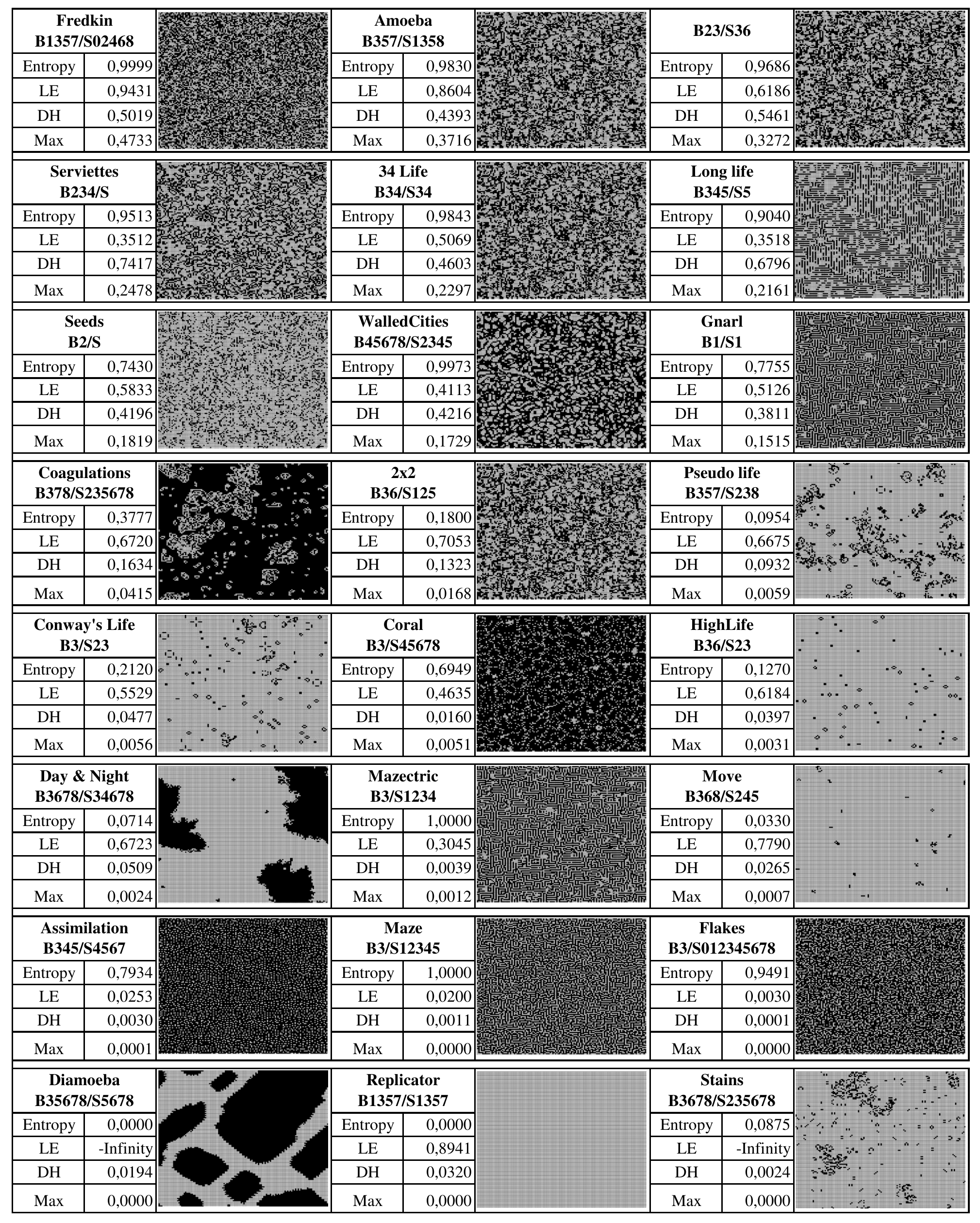}

\label{Tab:LifeLikeComplete}
\end{tabular}
\end{center}
\end{table*}

\subsection{CA as PRNG}
Most PRNGs are tested by suites or tools to analysis its behavior and truly randomness, thus according to Guan's experimental setup\cite{Guan}, all this experiments were tested with ENT\cite{ENT} and DIEHARD\cite{DIEHARD} suite.

The ENT test consist of five tests which are the Entropy test, Chi-square test, Serial correlation coefficient (SCC) test, Arithmetic Mean test, and Monte Carlo value for Pi, perhaps the complete details of the ENT tests is beyond the scope of this paper. In the other hand, the DIEHARD contains    18 different and independent statistical tests. The result of each test is called $p-value$. Most of the tests return a uniform p-value on $[0,1)$ and is considered a ``Fail" test, if the p-value is 1 or 0, in other cases can be considered as a ``Pass'' test \cite{DIEHARD}.

For both testes, was required a minimum of 10 Mb of random numbers, each PRNG received randomly chosen password and then executed until 10 Mb times. Were used 10 types of sets according to the number of blocks $\rho$, each one with 100 samples for testing.

In Table \ref{Tab:ENT10} is showed the average results with its standard deviation respectively for different values of $\rho$, all results are considered ``good'' PRNG because of the close approximation to theoretical randomness. Thus, for this proposal, we could take the $\rho=10$ blocks PRNG as the one with better results.

\begin{table*}[htbp]
\begin{center}
  \caption{The results of the ENT test suit, performed for different $\delta$ blocks}
\begin{tabular}{c}

\includegraphics[width=1\textwidth,height=0.35\textwidth]{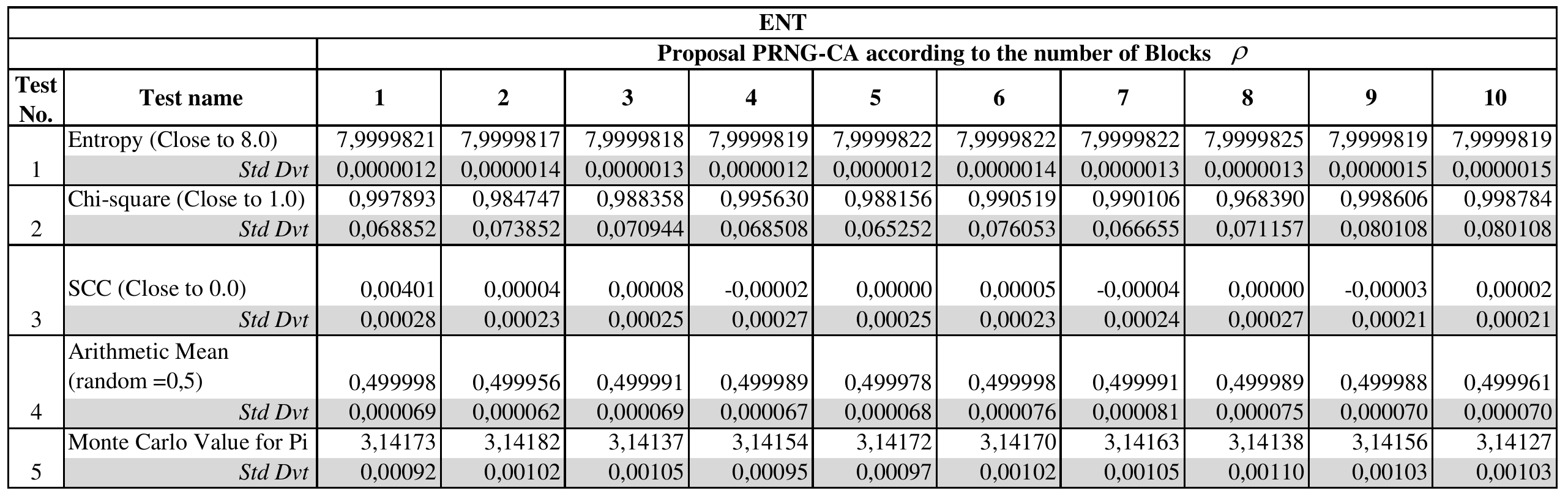}

\label{Tab:ENT10}
\end{tabular}
\end{center}
\end{table*}

In Table \ref{Tab:ENTComparison} is showed a comparison of the average results using the ENT suite, applied to our proposal based on $\rho=10$ blocks and others PRNG based on CA found in literature.
As it can be seen the proposed scheme is superior to other schemes, in Entropy test and the Serial Correlation Coefficient, although for the Chi-square does not seems to reaches to Shin et. al, \cite{Shin2009}.

\begin{table*}[htbp]
\begin{center}
  \caption{A comparison of the ENT test suite results}
\begin{tabular}{c}

\includegraphics[width=0.7\textwidth,height=0.15\textwidth]{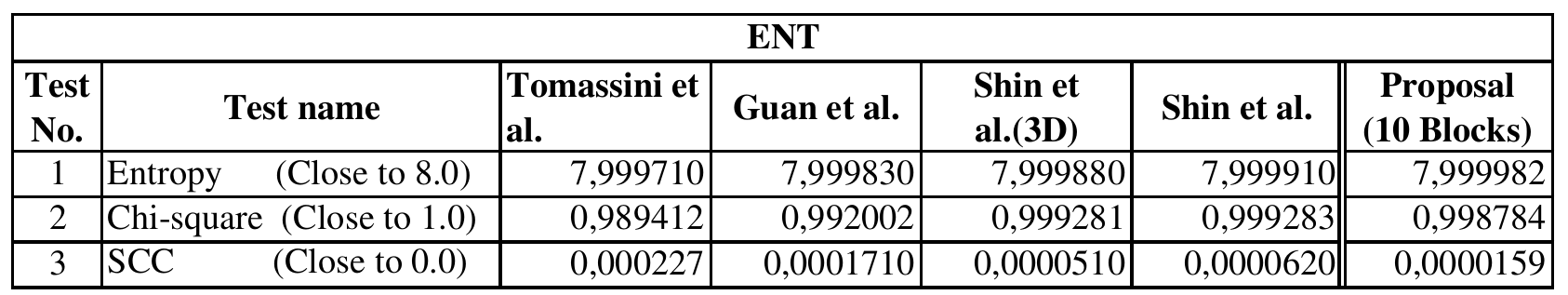}
\label{Tab:ENTComparison}
\end{tabular}
\\
Tomassini et. al \cite{TomassiniRandom} according to \cite{Shin2009};
\\
Guan et. al. (3D) \cite{Guan} according to \cite{Shin2009}; Shin et. al. (3D) \cite{Shin}; Shin et. al. \cite{Shin2009}
\end{center}
\end{table*}%


As with previous experiments, the experiments with different number of blocks on DIEHARD test, gave for each one a ``Pass'' result, thus it has been proven that the quality of randomness for $1,\ldots,10$ blocks got good results. So, in Table \ref{Tab:diehard} is showed a comparison of the DIEHARD test for our proposal by employing $\rho=10$ blocks and then compared with others PRNG based on CA found in literature.

Should be notice that the sequences of numbers generated by using just $4$ blocks, presents a sustain  number of blocks that improve cryptography thus by increasing the number of blocks will keeping this results, although will decrease the performance by slowing the encryption/decryption process.

\begin{table*}[htbp]
\begin{center}
  \caption{A comparison of the results of DIEHARD test in p-value pass rate 90\%}
\begin{tabular}{c}

\includegraphics[width=0.85\textwidth,height=0.5\textwidth]{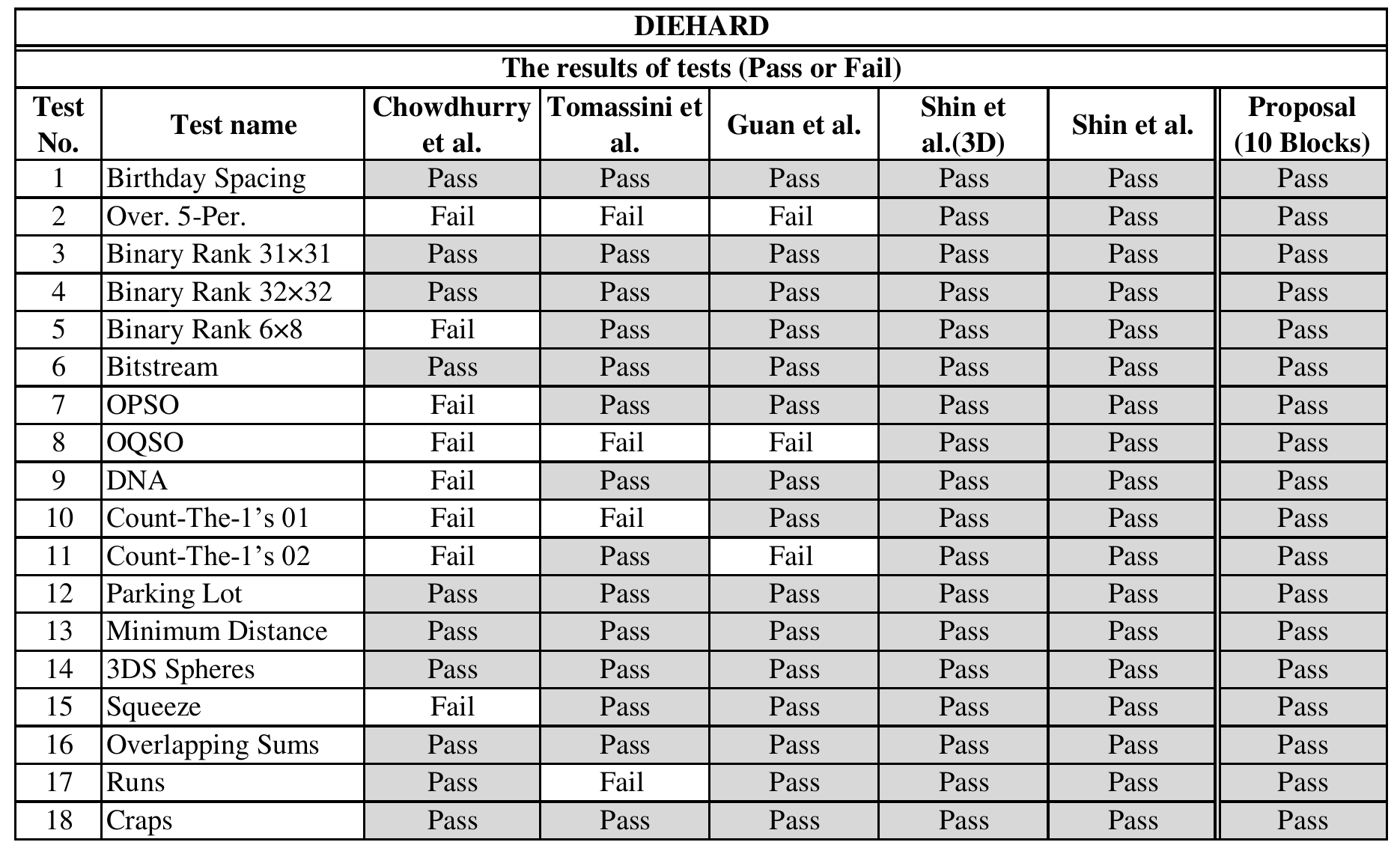}

\label{Tab:diehard}
\end{tabular}
\\Chowdhury et al. \cite{Chowdhuri} according to \cite{Shin2009}; Tomassini et. al \cite{TomassiniRandom} according to \cite{Shin2009};
\\
Guan et. al. (3D) \cite{Guan} according to \cite{Shin2009}; Shin et. al. (3D) \cite{Shin}; Shin et. al. \cite{Shin2009}
\end{center}
\end{table*}

\subsection{Encryption Results}
\begin{figure}[h]
\begin{center}$
\begin{array}{cc}
\includegraphics[width=0.2\textwidth,height=0.2\textwidth]{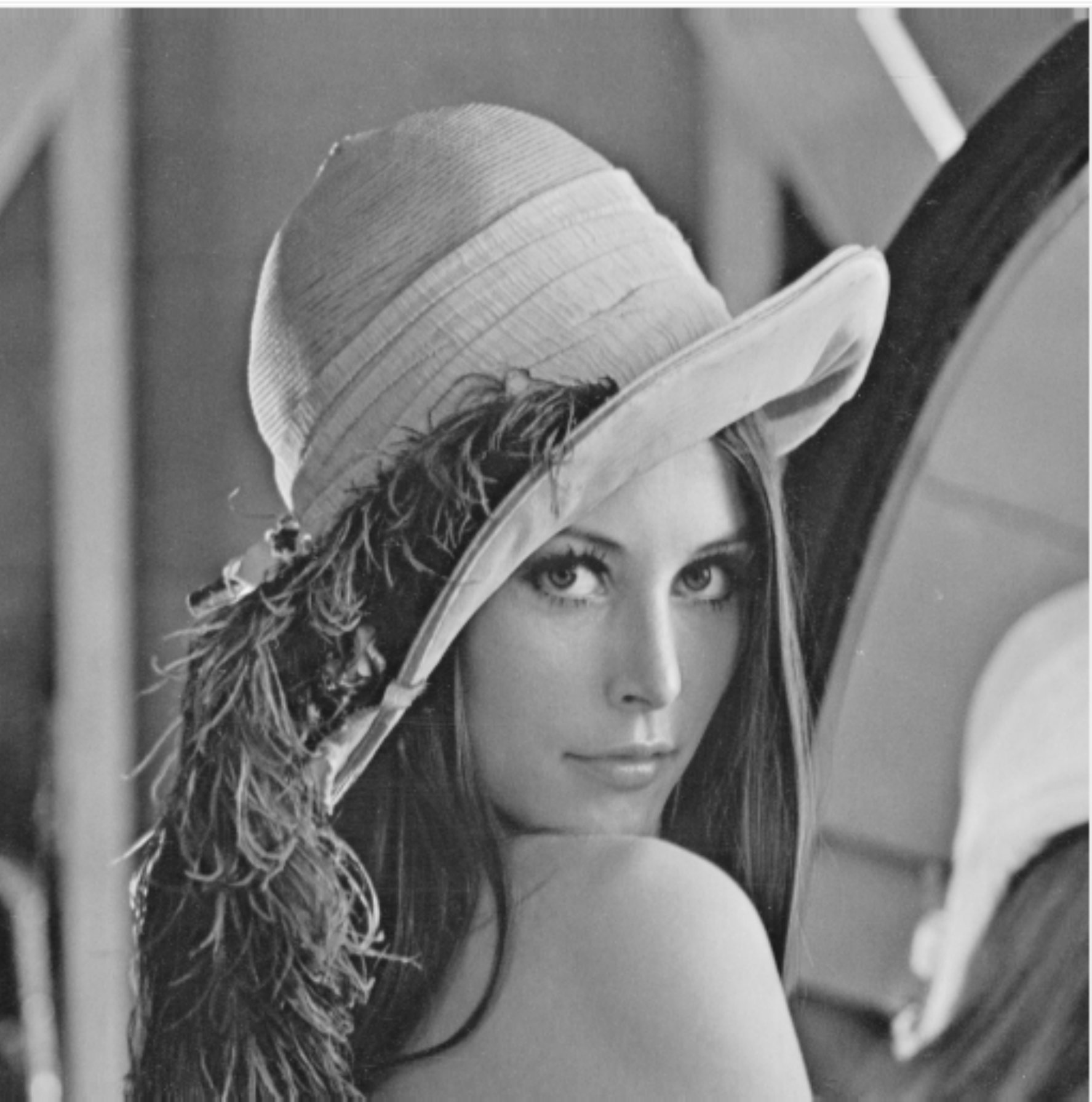}& \includegraphics[width=0.2\textwidth,height=0.2\textwidth]{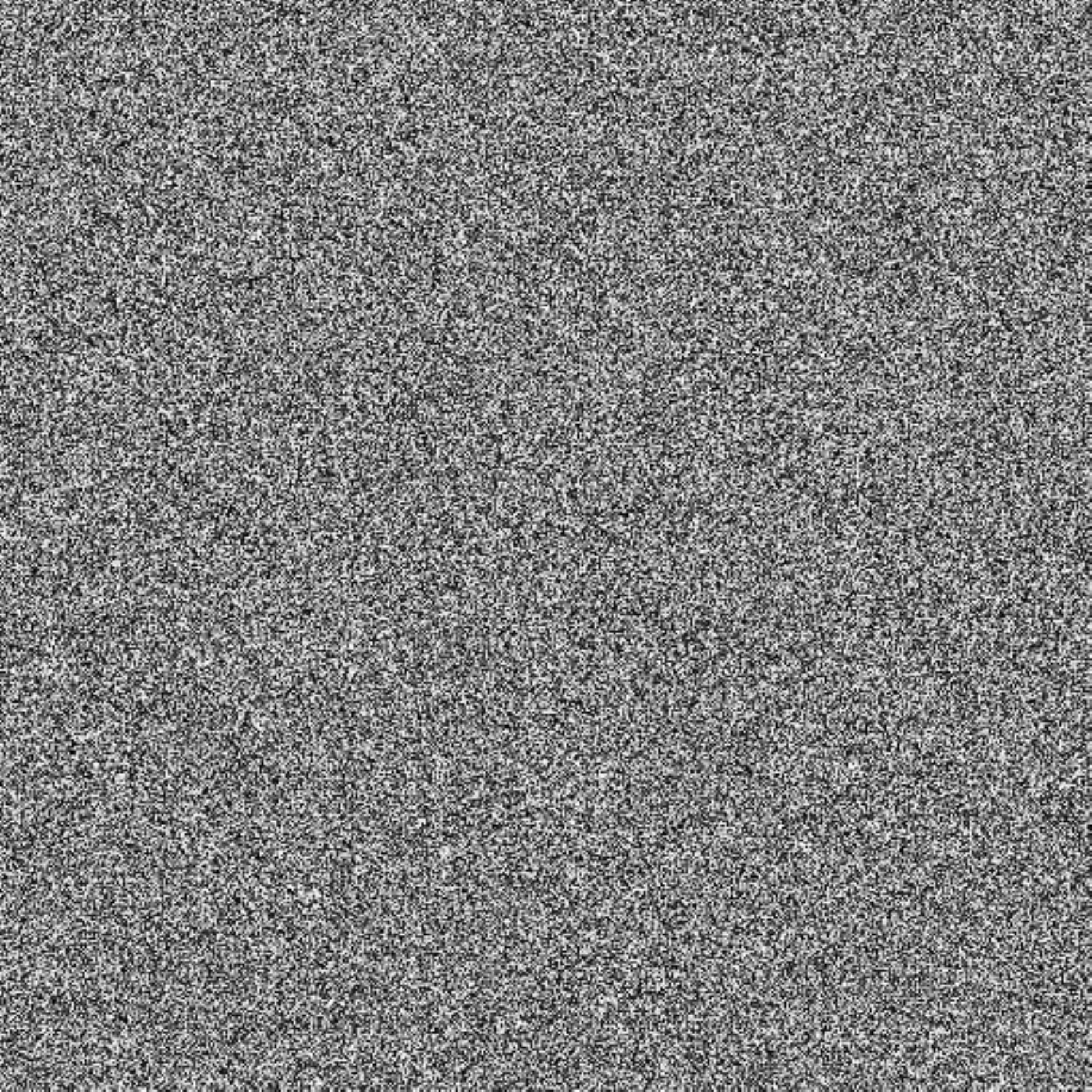}
\end{array}$
\end{center}
\caption{Left: The original plainimage. Right: The resultant cipherimage with the encryption method proposed }
\label{fig:plainimage}
\end{figure}

In this experiment the cipher was evaluated encrypting an image, on Figure \ref{fig:plainimage} is presented the original image (left), and the resultant cipher image (right). Visibly the ciphertext is  much more uniform and do not bears any resemblance of the plainimage.

\begin{figure}[h]
\begin{center}$
\begin{array}{cc}

\includegraphics[width=0.23\textwidth,height=0.18\textwidth]{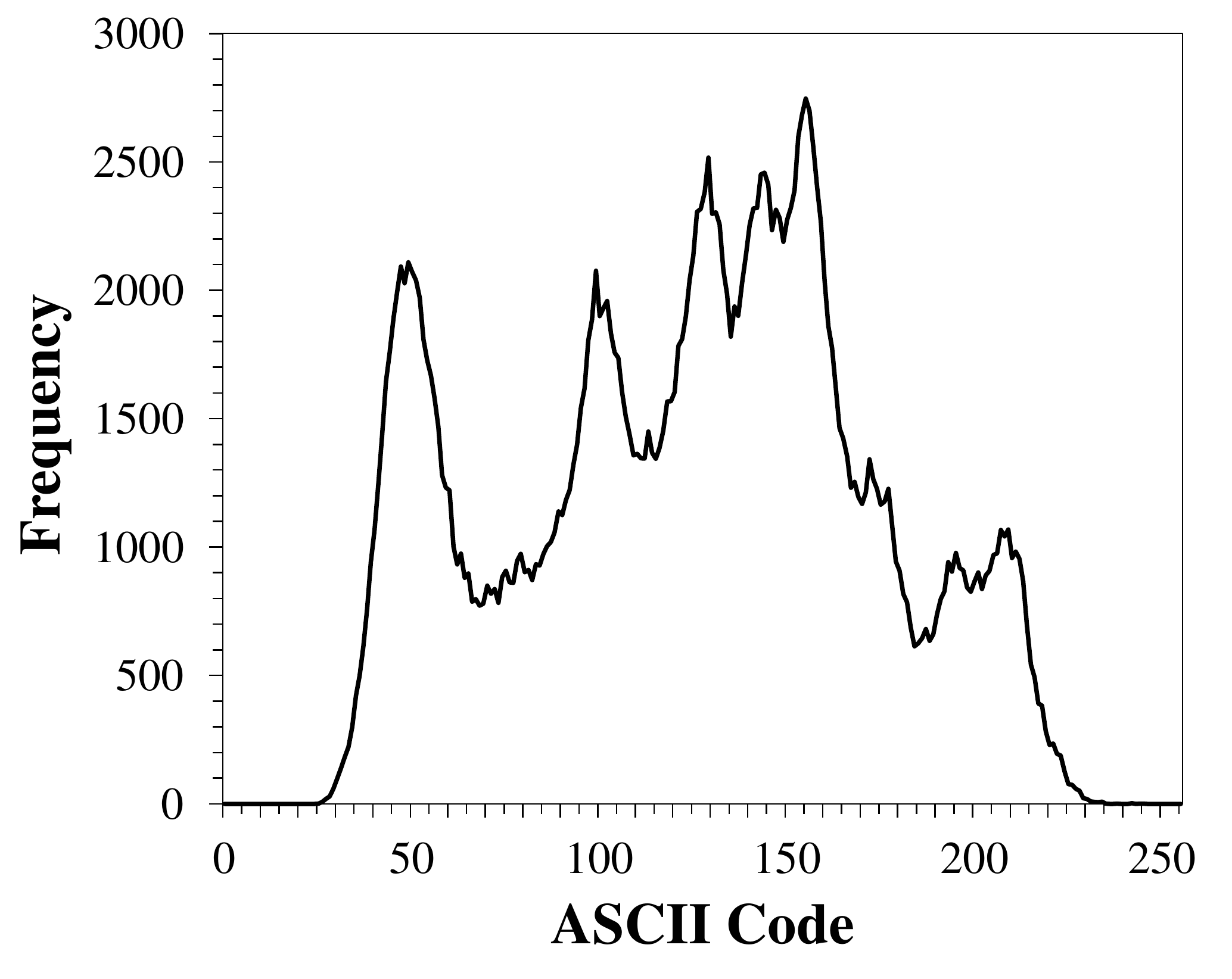}&
\includegraphics[width=0.23\textwidth,height=0.18\textwidth]{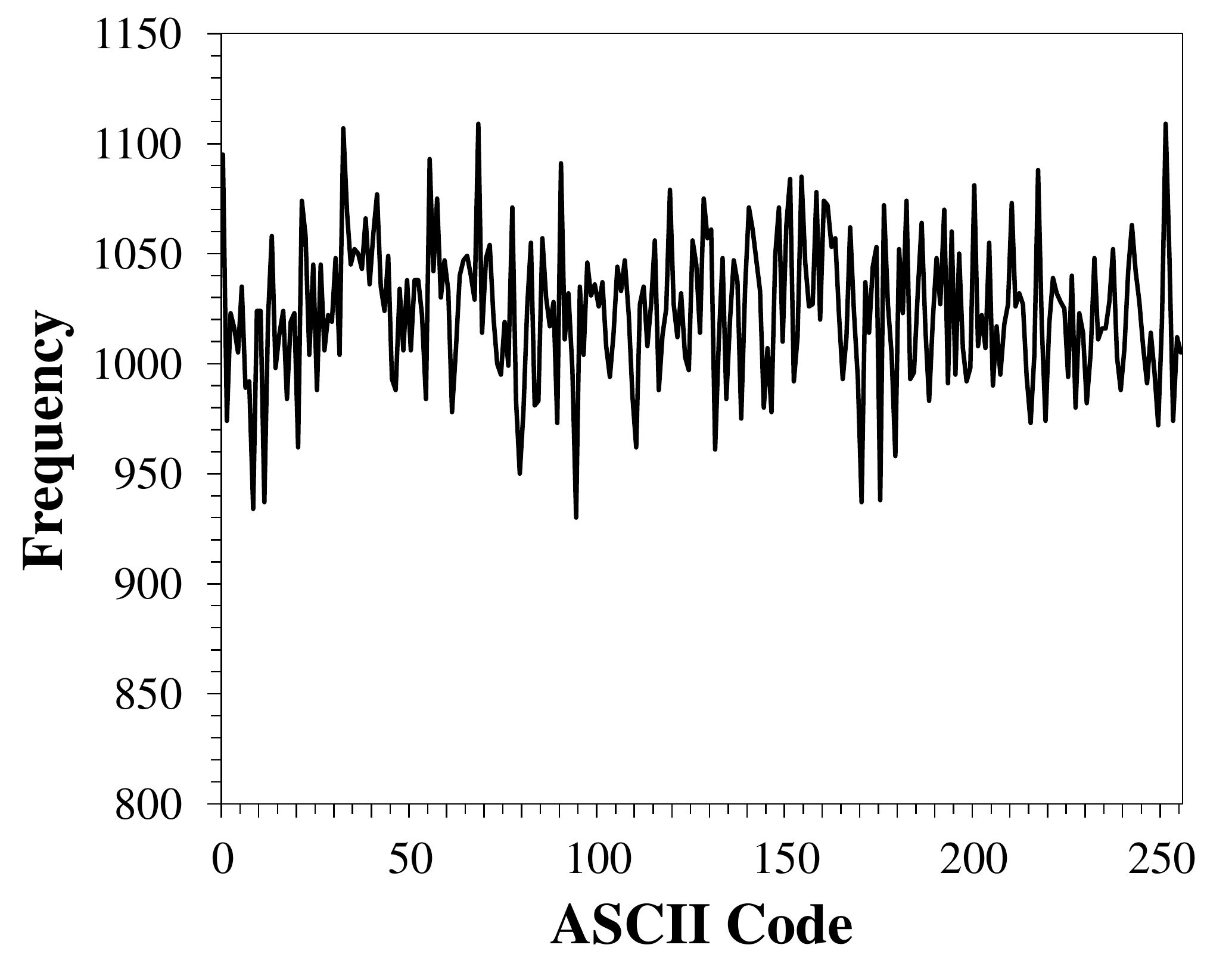}

\end{array}$
\end{center}
\caption{Frequency Analysis of the plainimage(left) and cipherimage(right) encrypted with the proposal method}
\label{fig:histogram}
\end{figure}

To analyze the statistical distribution of the PRNG, was generated an histogram of the plainimage and of the cipherimage which is shown on the left and right of Figure \ref{fig:histogram} respectively.
As it is shown the histogram of the cipherimage is much more regular and almost plain, which demonstrated that all the 256 numbers(ASCCI code) generated with the CA are scattered all over, with almost same percentages regions, by means, hiding information.

\subsection{Fourier's Power Spectrum}

\begin{figure}[h]
\begin{center}$
\begin{array}{cc}
\includegraphics[width=0.2\textwidth,height=0.2\textwidth]{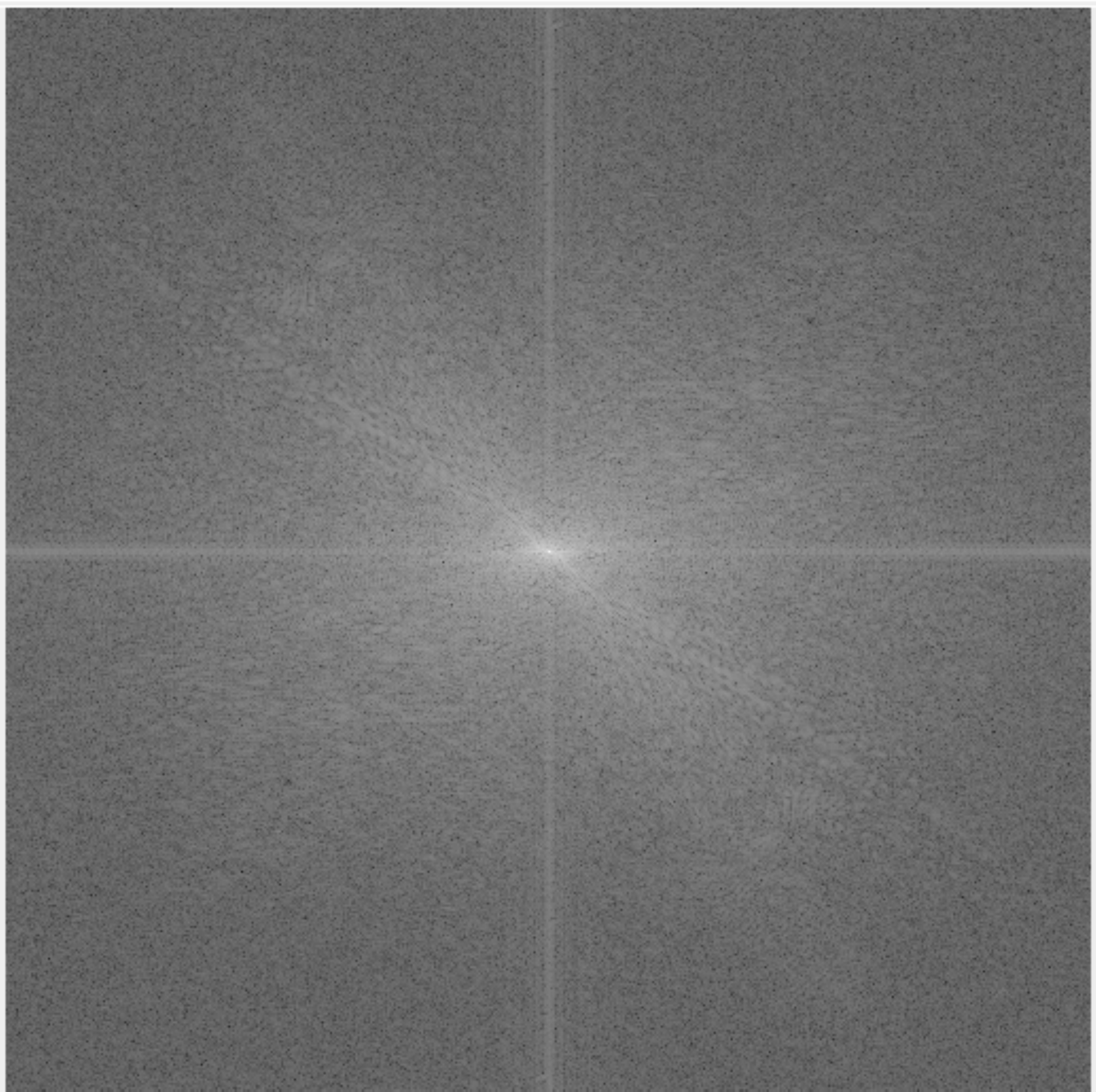}&
\includegraphics[width=0.2\textwidth,height=0.2\textwidth]{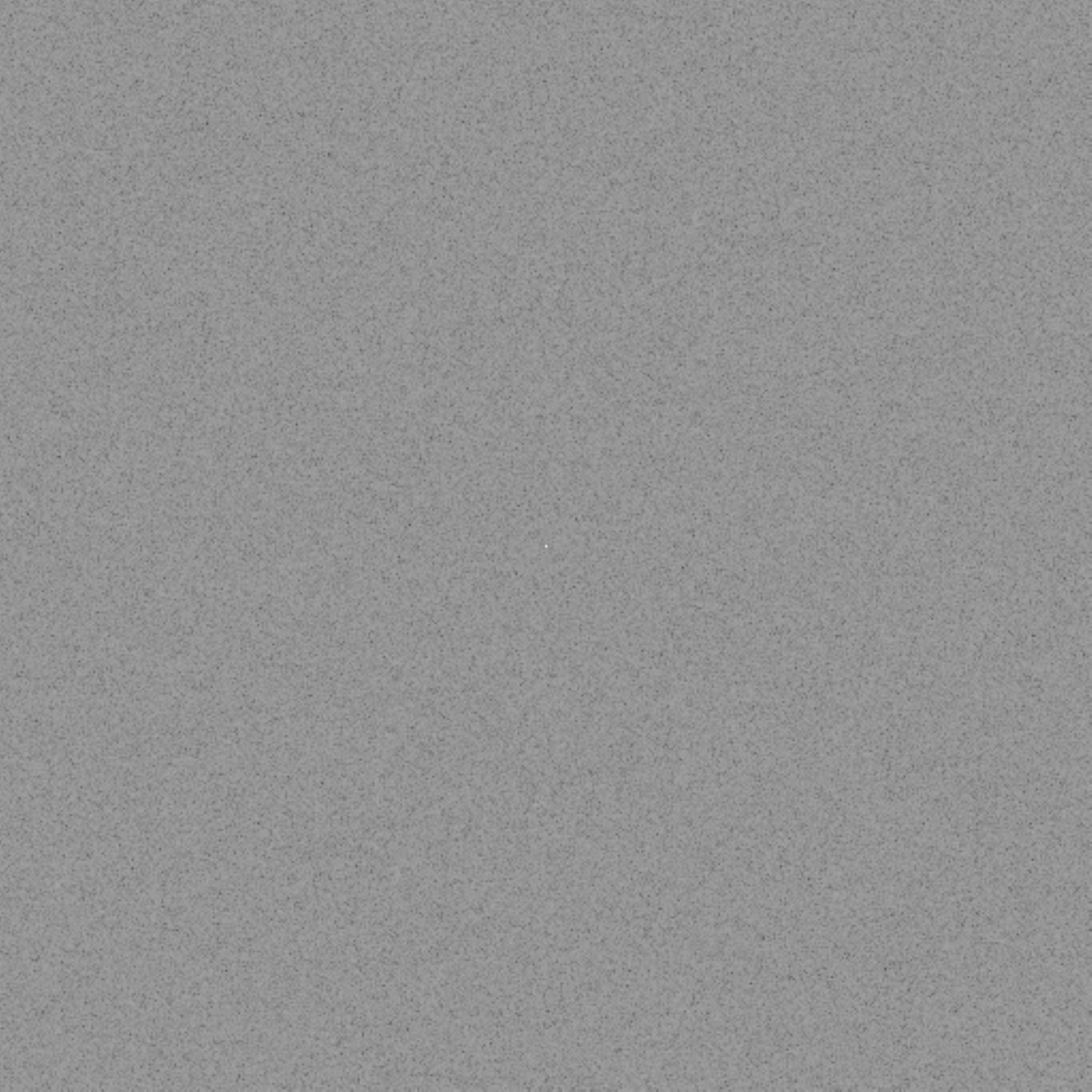}
\end{array}$
\end{center}
\caption{2D Fourier power spectrum for the plainimage(left) and cipherimage(right)}
\label{Fourier}
\end{figure}

On Figure \ref{Fourier}, for both the plainimage and cipherimage, has been performed and compared to 2D Fourier power spectrum, as we can seen the spectrum of the cipherimage minimal presents a white noise, which it is said to presents any frequency
information,  demonstrating any compressible information can be achieved in the
cipher image\cite{Anderson}.

\section{Conclusions}
In this paper, a chaotic encryption method based on the ``Life-Like'' Cellular Automata was proposed. The cryptosystem we have described is based on two-dimensional CA which is able to generate pseudo-random numbers of very good quality as measured by both: chaotic measures(Lyapunov Exponent, Entropy, Hamming distance, and combination of them) and statistical tests(DIEHARD and ENT test suites), as they achieve better randomness quality than others based on CA compared in literature.

Here were described  chaotic properties of CA as good reasons to be employed in cryptography, furthermore, its simplicity and low cost of implementation as one of its main advantages.


\section*{Acknowledgements}
M. J. M. acknowledges support from FAPESP (The State of São Paulo Research Foundation) (2011/05461-0).
A. G. M. acknowledges the Brazilian agency CNPq (XXXX) for support. 
O.M.B. acknowledges support from CNPq (Grant \#308449/2010-0 and \#473893/2010-0) and FAPESP (Grant \# 2011/01523-1). 



\begin{thebibliography}{10}
\providecommand{\url}[1]{#1}
\csname url@samestyle\endcsname
\providecommand{\newblock}{\relax}
\providecommand{\bibinfo}[2]{#2}
\providecommand{\BIBentrySTDinterwordspacing}{\spaceskip=0pt\relax}
\providecommand{\BIBentryALTinterwordstretchfactor}{4}
\providecommand{\BIBentryALTinterwordspacing}{\spaceskip=\fontdimen2\font plus
\BIBentryALTinterwordstretchfactor\fontdimen3\font minus
  \fontdimen4\font\relax}
\providecommand{\BIBforeignlanguage}[2]{{%
\expandafter\ifx\csname l@#1\endcsname\relax
\typeout{** WARNING: IEEEtran.bst: No hyphenation pattern has been}%
\typeout{** loaded for the language `#1'. Using the pattern for}%
\typeout{** the default language instead.}%
\else
\language=\csname l@#1\endcsname
\fi
#2}}
\providecommand{\BIBdecl}{\relax}
\BIBdecl

\bibitem{CrackDES}
F.~E. Frontier, \emph{Cracking DES: Secrets of Encryption Research, Wiretap
  Politics and Chip Design}, M.~Loukides and J.~Gilmore, Eds.\hskip 1em plus
  0.5em minus 0.4em\relax O'Reilly \& Associates, Inc., 1998.

\bibitem{CrackMD5}
\BIBentryALTinterwordspacing
T.~Xie, D.~Feng, and F.~Liu, ``A new collision differential for md5 with its
  full differential path,'' IACR Cryptology ePrint Archive, 2008. [Online].
  Available: \url{http://eprint.iacr.org/2008/230}
\BIBentrySTDinterwordspacing

\bibitem{CrackSHA1}
M.~St\'{e}phane, ``Classification and generation of disturbance vectors for
  collision attacks against sha-1,'' \emph{IJDCC Des. Codes Cryptography},
  vol.~59, pp. 247--263, 2011.

\bibitem{neumann}
J.~V. Neumann, \emph{Theory of Self-Reproducing Automata}, A.~W. Burks,
  Ed.\hskip 1em plus 0.5em minus 0.4em\relax University of Illinois Press,
  1966.

\bibitem{Wolfram}
S.~Wolfram, \emph{A New Kind of Science}.\hskip 1em plus 0.5em minus
  0.4em\relax Wolfram Media, Inc, 2002.

\bibitem{Xiang}
X.~T., L.~X., T.~G., C.~Y., and W.~K.W., ``A novel block cryptosystem based on
  iterating a chaotic map,'' \emph{Phys. Lett. A.}, vol. 349, pp. 109--115,
  2006.

\bibitem{Anderson}
A.~Marco, A.~Martinez, and O.~Bruno, ``Fast, parallel and secure cryptography
  algorithm using loren'z attractor,'' \emph{IJMPC Mod Phy C}, vol.~21, pp.
  265--382, 2010.

\bibitem{Deng}
S.~Deng, D.~Xiao, Y.~Li, and W.~Peng, ``A novel combined cryptographic and hash
  algorithm based on chaotic control character,'' \emph{Communications in
  Nonlinear Science and Numerical Simulation}, vol.~14, no.~11, pp. 3889--3900,
  2009.

\bibitem{Alvarez2006}
G.~Alvarez and S.~Li, ``Some basic cryptographic requirements for chaos-based
  cryptosystems,'' \emph{IJBC}, vol.~16, pp. 2129--2151, 2006.

\bibitem{Habutsu}
T.~Habutsu, Y.~Nishio, I.~Sasase, and S.~Mori, ``A secret key cryptosystem by
  iterating a chaotic map,'' in \emph{Theory and Application of Cryptographic
  Techniques}, 1991, pp. 127--140.

\bibitem{Nandi}
S.~Nandi, B.~K. Kar, and P.~P. Chaudhuri, ``Theory and applications of cellular
  automata in cryptography,'' \emph{IEEE Transactions on Computers}, vol.~43,
  pp. 1346--1357, 1994.

\bibitem{Tomassini}
M.~Tomassini and M.~Perrenoud, ``Stream cyphers with one and two-dimensional
  cellular automata,'' in \emph{Proceedings of the 6th International Conference
  on Parallel Problem Solving from Nature}.\hskip 1em plus 0.5em minus
  0.4em\relax Springer-Verlag, 2000, pp. 722--731.

\bibitem{Seredynski}
F.~Seredynski, P.~Bouvry, and A.~Y. Zomaya, ``Cellular automata computations
  and secret key cryptography,'' \emph{Parallel Computing}, vol.~30, pp.
  753--766, 2004.

\bibitem{TomassiniRandom}
M.~Tomassini, M.~Sipper, and M.~Perrenoud, ``On the generation of high-quality
  random numbers by two-dimensional cellular automata,'' \emph{Computers, IEEE
  Transactions on}, vol.~49, no.~10, pp. 1146--1151, 2000.

\bibitem{Martin}
B.~Martin and P.~Solé, ``Pseudorandom sequences generated by cellular
  automata,'' \emph{Computing Research Repository}, vol. abs/0807.3, 2008.

\bibitem{Guan}
S.~U. Guan and S.~K. Tan, ``Pseudorandom number generation with
  self-programmable cellular automata,'' \emph{IEEE Transactions on
  Computer-aided Design of Integrated Circuits and Systems}, vol.~23, pp.
  1095--1101, 2004.

\bibitem{Shin}
S.~ho~Shin, G.~dal Park, and K.~young Yoo, ``A virtual three-dimension cellular
  automata pseudorandom number generator based on the moore neighborhood
  method,'' in \emph{International Conference on Intelligent Computing}.\hskip
  1em plus 0.5em minus 0.4em\relax Springer-Verlag, 2008, pp. 174--181.

\bibitem{Shin2009}
S.-H. Shin and K.-Y. Yoo, ``An efficient prng based on the hybrid between one-
  and two-dimensional cellular automata,'' in \emph{ITNG'09}, 2009, pp.
  498--503.

\bibitem{Seth}
A.~Seth, S.~Bandyopadhyay, and U.~Maulik, ``Pseudorandom pattern generation by
  a 4neighborhood cellular automata based on a probabilistic analysis,'' in
  \emph{Proceedings of the International MultiConference of Engineers and
  Computer Scientists IMECS}, vol.~2, 2008.

\bibitem{Wolf}
A.~Wolf, J.~B. Swift, H.~L. Swinney, and J.~A. Vastano, ``Determining lyapunov
  exponents from a time series,'' \emph{Physica}, pp. 285--317, 1985.

\bibitem{DIEHARD}
G.~Marsaglia, ``The marsaglia random number cdrom, with the diehard battery of
  tests of randomness,'' Florida, USA, 1998.

\bibitem{ENT}
J.~Walker, ``A pseudorandom number sequence test program,'' 1998.

\bibitem{Tanenbaum}
A.~S. Tanenbaum, \emph{Computer Networks}, 4th~ed.\hskip 1em plus 0.5em minus
  0.4em\relax Prentice Hall Professional Technical Reference, 2003.

\bibitem{Stallings}
W.~Stallings, \emph{Cryptography and Network Security : Principles and
  Practice}, 5th~ed.\hskip 1em plus 0.5em minus 0.4em\relax Prentice Hall,
  2002.

\bibitem{Baptista}
M.~Baptista, ``Cryptography with chaos,'' \emph{Phys Lett A}, vol. 240, pp.
  50--54, 1998.

\bibitem{Gutowitz}
H.~Gutowitz, ``Cryptography with dynamical systems,'' in \emph{NATO Advanced
  Study Institute}, 1996.

\bibitem{Golly}
A.~Trevorrow and T.~Rokicki, ``Golly an open source, cross-platform application
  for exploring conway's game of life and other cellular automata.''

\bibitem{Adamatzky}
A.~Adamatzky, \emph{Game of Life Cellular Automata}, 1st~ed.\hskip 1em plus
  0.5em minus 0.4em\relax Springer Publishing Company, Incorporated, 2010.

\bibitem{BookChaos}
M.~Schroeder, \emph{Fractals, Chaos, Power Laws: : Minutes from an Infinite
  Paradise}, 1st~ed., W.~H. Freeman and Company, Eds.\hskip 1em plus 0.5em
  minus 0.4em\relax Freeman, 1991.

\bibitem{Guckenheimer}
J.~Guckenheimer and P.~Holmes, \emph{Nonlinear Oscillations, Dynamical Systems,
  and Bifurcations of Vector Fields}, 1st~ed.\hskip 1em plus 0.5em minus
  0.4em\relax Springer-Verlag, 1983, vol.~42.

\bibitem{EckmannRuelle}
J.~P. Eckmann and D.~Ruelle, ``Ergodic theory of chaos and strange
  attractors,'' \emph{Rev. Mod. Phys.}, vol.~57, no. 3, part 1, pp. 617--656,
  1985.

\bibitem{Tisseur2003}
P.~Tisseur, ``Cellular automata and lyapunov exponents,'' 2003,
  http://arXiv:math/0312136v1.

\bibitem{Tisseur2005}
P.~Tisseur, ``Always finite entropy and lyapunov exponents of two-dimensional
  cellular automata,'' 2005, http://arXiv:math/0502440v1.

\bibitem{Cattaneo}
G.~Cattaneo, E.~Formenti, L.~Margara, and G.~Mauri, ``On the dynamical behavior
  of chaotic cellular automata,'' \emph{Theoretical Computer Science}, vol.
  217, pp. 31--51, 1999.

\bibitem{XUXU}
X.~Xu, Y.~Song, and S.~P. Banks, ``On the dynamical behavior of cellular
  automata,'' \emph{IJBC}, vol.~19, 2009.

\bibitem{Shimada}
I.~Shimada and T.~Nagashima, ``A numerical approach to ergodic problem of
  dissipative dynamical systems,'' \emph{Progress of Theoretical Physics},
  vol.~61, pp. 1605--1616, 1979.

\bibitem{Benettin}
G.~Benettin, L.~Galgani, A.~Giorgilli, and J.-M. Strelcyn, ``Lyapunov
  characteristic exponents for smooth dynamical systems and for hamiltonian
  systems; a method for computing all of them. part 1: Theory,''
  \emph{Meccanica}, vol.~15, pp. 9--20, 1980.

\bibitem{Baetens2010}
J.~M. Baetens and B.~De~Baets, ``Phenomenological study of irregular cellular
  automata based on lyapunov exponents and jacobians.'' \emph{Chaos}, vol.~20,
  no.~3, p. 033112, 2010.

\bibitem{Chowdhuri}
D.~R. Chowdhuri, P.~Subbarao, and P.~P. Chaudhuri, ``Characterization of
  two-dimensional cellular automata using matrix algebra,'' \emph{Inf. Sci.},
  vol.~71, pp. 289--314, 1993.

\end{thebibliography}


\end{document}